\documentclass[reqno]{amsart}
\usepackage{setspace}
\usepackage{graphicx}
\usepackage{amsmath}
\usepackage{amsfonts}
\usepackage{amssymb}
\begin{document}

\centerline{\huge \bf A Computer Program for Borsuk's Conjecture}

\bigskip\medskip
\centerline{\Large\bf Chuanming Zong}

\vspace{1cm}
\centerline{\begin{minipage}{12.8cm}
{\bf Abstract.} In 1933, Borsuk proposed the following problem: Can every bounded set in $\mathbb{E}^n$ be divided into $n+1$ subsets of smaller diameters?  This problem has been studied by many authors, and a lot of partial results have been discovered. In particular, Kahn and Kalai's counterexamples surprised the mathematical community in 1993. Nevertheless, the problem is still far away from being completely resolved. This paper presents a broad review on related subjects and, based on a novel reformulation, introduces a computer proof program to deal with this well-known problem.
\end{minipage}}


\vspace{1cm}
\noindent
{\LARGE\bf 1. Borsuk's Conjecture}

\bigskip
\noindent
Let $X$ be a subset of the $n$-dimensional Euclidean space $\mathbb{E}^n$ with diameter
$$d(X)=\sup_{{\bf x},\, {\bf y}\in X} \| {\bf x}, {\bf y}\|,\eqno(1.1)$$
where $\| {\bf x}, {\bf y}\|$ denotes the Euclidean distance between ${\bf x}$ and ${\bf y}$. As usual, let $\partial (X)$ and ${\rm int}(X)$ denote the boundary and the interior of $X$, respectively.

\medskip
In 1933, K. Borsuk \cite{bors33} studied the continuous maps between metric spaces. As a corollary of his main result it was shown that, {\it whenever an $n$-dimensional Euclidean ball is divided into $n$ subsets, at least one of these subsets has the same diameter of the ball.} Then, at the end of the paper, he proposed the following problem:

\medskip
\noindent
{\bf Borsuk's Problem.} {\it Can every bounded set in $\mathbb{E}^n$ be divided into $n+1$ subsets of smaller diameter?}

\medskip
Usually, the positive statement of this problem is referred as {\it Borsuk's Conjecture}, though Borsuk himself only proposed it as an open problem.
For convenience, let $b(X)$ denote the smallest number such that $X$ can be partitioned into $b(X)$ subsets of smaller diameter. Then, Borsuk's conjecture can be reformulated as following.

\medskip\noindent
{\bf Borsuk's Conjecture.} {\it For every bounded subset $X$ of the $n$-dimensional Euclidean space $\mathbb{E}^n$ we have}
$$b(X)\le n+1.$$

\vspace{1cm}
\noindent
{\LARGE\bf 2. Reductions and Positive Results}

\bigskip
\noindent
{\bf Definition 1.} Assume that $X$ is a subset of $\mathbb{E}^n$. Then we define
$$\widehat{X}=\left\{\lambda {\bf x}_1+(1-\lambda ){\bf x}_2: {\bf x}_i\in K,\ 0\le \lambda \le 1\right\}.$$
Usually, we call $\widehat{X}$ the convex hull of $X$. In particular, we call a compact subset $K$ of $\mathbb{E}^n$ an $n$-dimensional convex body if it has nonempty interior and satisfying $K=\widehat{K}$.

\medskip
It is obvious that $X\subseteq \widehat{X}$ and $d(X)=d(\widehat{X})$. Therefore, to solve Borsuk's problem, it is sufficient to deal with all the convex bodies $K$.

\medskip\noindent
{\bf Definition 2.} Assume that $C$ is an $n$-dimensional convex body and ${\bf u}$ is a unit vector in $\mathbb{E}^n$. Then $C$ has two tangent hyperplanes $H_1$ and $H_2$ with norm ${\bf u}$. Let $d(C,{\bf u})$ denote the distance between $H_1$ and $H_2$. If there is a constant $c$ such that
$$d(C,{\bf u})=c$$
holds for all unit vectors ${\bf u}$, we will call $C$ a convex body of constant width.

\medskip
Clearly, balls are convex bodies of constant width. In addition, Reuleaux triangles and Meissner bodies are particular examples. In fact, convex bodies of constant width in $\mathbb{E}^n$ is a fascinating field of research. There are hundreds of papers on this subject. Many well-known mathematicians such as W. Blaschke, M. Fujiwara, H. Lebesgue, K. Reidemeister, L. A. Santal\'o and W. S\"uss have made contribution to this field. Nevertheless, up to now some basic problems about convex bodies of constant width are still open (see \cite{chak83, mart19}). The next result is useful for Borsuk's problem.

\medskip\noindent
{\bf Lemma 1 (P$\acute{\bf a}$l \cite{pal20}, Lebesgue \cite{lebe21}).} {\it For every bounded set $X$ in $\mathbb{E}^n$ there is a convex body $C$ of constant width satisfying both $X\subseteq C$ and}
$$d(X)=d(C).$$

\begin{figure}[ht]
\includegraphics[scale=0.5]{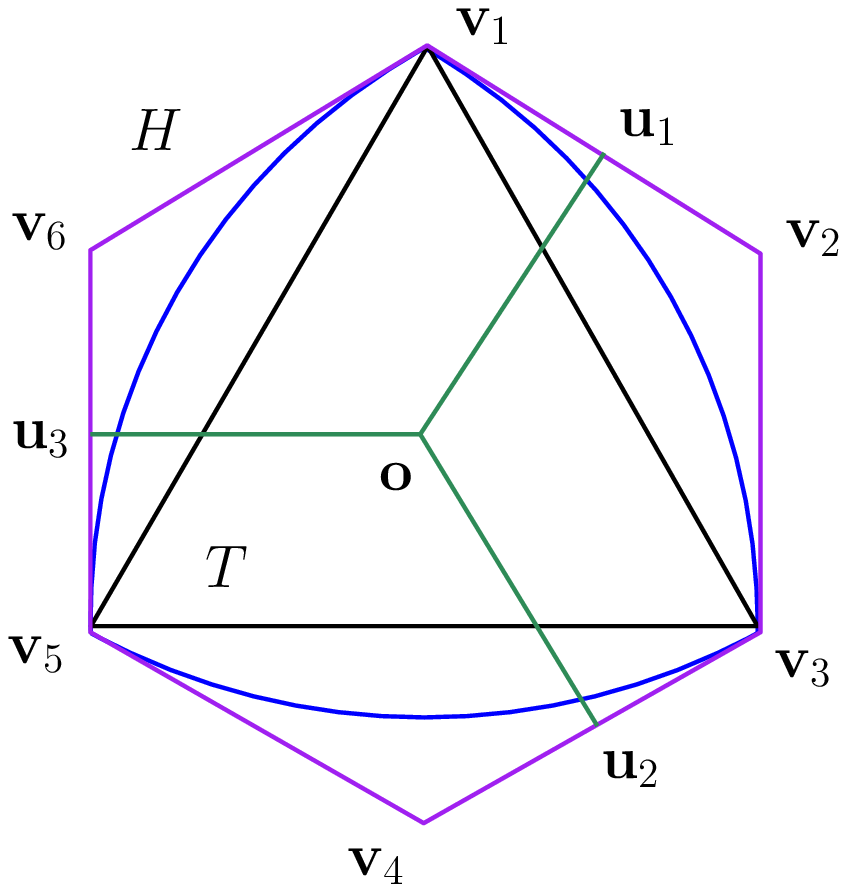}
\end{figure}

\centerline{\begin{minipage}{11cm}
{\bf Figure 1.} {\it A regular triangle can be embedded into a Reuleaux triangle, the Reuleaux triangle can be embedded into a regular hexagon, and the regular hexagon can be divided into three subsets of smaller diameter.}
\end{minipage}}

\bigskip
Based on this lemma, to prove Borsuk's conjecture it is sufficient to deal with all the convex bodies of unit constant width. On the other hand, by continuity argument, one can deduce that every convex body of unit constant width can be inscribed into a regular hexagon that the distance between the opposite sides is $1$. In other words, the hexagon has edge length $2/\sqrt{3}$. Then, it is easy to see that the hexagon can be divided into three parts of diameter ${{\sqrt{3}}/ 2}$, as shown by Figure 1. Thus one obtains the following theorem.

\medskip\noindent
{\bf Theorem 1 (Bonnesen and Fenchel \cite{bonn34}).} {\it Every two-dimensional set of diameter $d$ can be divided into three subsets of diameter at most ${{\sqrt{3}}\over 2}d$.}

\medskip\noindent
{\bf Remark 1.} Clearly, the constant ${{\sqrt{3}}\over 2}\approx 0.866\ldots $ is optimal.

\medskip
In 1945, H. Hadwiger \cite{hadw45} claimed a proof for Borsuk's conjecture based on Lemma 1. Soon afterwards, he realized that his proof was relied on the assumption that convex bodies of constant width have regular boundaries, which is apparently wrong. In fact, he was able to prove the following result.

\medskip\noindent
{\bf Theorem 2 (Hadwiger \cite{hadw46}).} {\it Every $n$-dimensional convex body with smooth boundary can be divided into $n+1$ subsets of smaller diameter.}

\medskip
Let $B$ denote the $n$-dimensional unit ball centered at the origin ${\bf o}$ and assume that $K$ is an $n$-dimensional convex body with a smooth boundary and ${\bf o}\in {\rm int}(K)$. Let ${\bf x}$ be a boundary point of $K$, let ${\bf u}({\bf x})$ denote the unit norm of $K$ at ${\bf x}$, and define
$$f({\bf x})={\bf u}({\bf x})-{\bf x}\eqno(2.1)$$
to be a map from $\partial (K)$ to $\partial (B)$. Assume that $Y_1$, $Y_2$, $\ldots $, $Y_{n+1}$ are $n+1$ subsets of $\partial (B)$ satisfying both
$$\partial (B)=\bigcup_{i=1}^{n+1}Y_i\eqno(2.2)$$
and
$$d(Y_i)<2, \quad i=1, 2, \ldots, n+1.\eqno(2.3)$$
Let $X_1$, $X_2$, $\ldots $, $X_{n+1}$ be subsets of $\partial (K)$ satisfying
$$f(X_i)=Y_i\eqno(2.4)$$
and define $Z_i$ to be the convex hull of ${\bf o}\cup X_i$. Then, it can be easily shown that
$$d(Z_i)<d(K)\eqno(2.5)$$
holds for all $i=1$, $2$, $\ldots$, $n+1$. Hadwiger's theorem is proved.

\medskip
In 1947, J. Perkal \cite{perk47} claimed that Borsuk's conjecture was correct in $\mathbb{E}^3$. However, he did not give the proof idea. In 1955, by modifying Hadwiger's approach, H. G. Eggleston \cite{eggl55} presented a detailed proof for the three-dimensional case of the conjecture. Later, different proofs were discovered by B. Gr\"unbaum \cite{grun57}, A. Happes and P. R\'ev\'esz \cite{hepp56}, and A. Happes \cite{hepp57}.

\medskip\noindent
{\bf Theorem 3 (Perkal \cite{perk47}, Eggleston \cite{eggl55}).} {\it Every three-dimensional bounded set can be divided into four subsets of smaller diameters.}

\medskip
In fact, the reduction idea can be extend to three dimensions (see Gr\"unbaum \cite{grun57}). First, every set of diameter one is a subset of a convex body of diameter one. Second, every convex body of diameter one is a subset of a set of constant width one. Third, every set of unit constant width can be embedded into a regular octahedron whose opposite facets are distance one apart. Fourth, by cutting off three suitable small pyramids from the octahedron one obtains a suitable truncated octahedron, as shown by Figure 2. Finally, the truncated octahedron can be divided into four polytopes of diameter less than $0.9887.$ Thus, Theorem 3 is proved.

\begin{figure}[ht]
\centering
\includegraphics[height=5cm,width=5cm,angle=0]{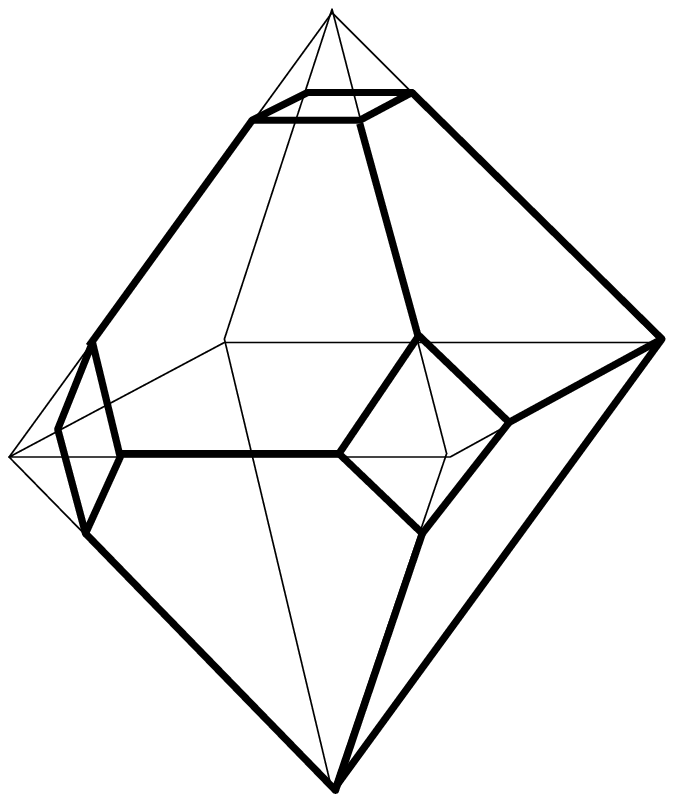}
\end{figure}

\centerline{\begin{minipage}{11.3cm}
{\bf Figure 2.} {\it The truncated octahedron, which contains the set of unit constant width, can be divided into four subsets of diameter smaller than one.}
\end{minipage}}

\bigskip\noindent
{\bf Remark 2.} It is easy to see that the constant $0.9887$ is not optimal. Clearly, to determine the optimal constant is a challenging problem.

\medskip
In attacking Borsuk's problem in higher dimensions, many partial results have been achieved. In 1945, H. Hadwiger \cite{hadw45} found that $b(K)\le n+1$ for every n-dimensional smooth convex body $K$ (Theorem 2). In 1955, H. Lenz \cite{lenz55} showed that, on the one hand, $b(K)\le n$ if $K$
has a smooth boundary but nonconstant width, while on the other hand, $b(K)\ge n+1$ for all sets of constant width. In 1971, A. S. Riesling \cite{ries71} showed  that $b(K)\le n+1$ for every n-dimensional centrally symmetric  convex body $K$. In 1971, C. A. Rogers \cite{roge71} proved that $b(K)\le n+1$ when $K$ is invariant under the {\it symmetry group} of a regular n-dimensional simplex.

\vspace{1cm}
\noindent
{\LARGE\bf 3. Counterexamples}

\bigskip
\noindent
In 1993, J. Kahn and G. Kalai \cite{kahn93} made the following counterintuitive discovery.

\medskip\noindent
{\bf Theorem 4 (Kahn and Kalai \cite{kahn93}).} {\it For every integer $n\ge 1$, there exists a subset $X_n$ of the $n$-dimensional Euclidean space such that}
$$b(X_n)\ge 1.07^{\sqrt{n}}.$$

\medskip\noindent
{\bf Remark 3.} Clearly, $1.07^{\sqrt{n}}$ is much larger than $n+1$ when $n$ is sufficiently large. In particular, whenever $n>21800$, we have
$$1.07^{\sqrt{n}}>n+1,$$
which gives counterexamples to Borsuk's conjecture in high dimensions.

\bigskip
Kahn and Kalai's counterexamples were indeed surprising. However, the proof idea is very natural, once understood. In 1981, D. Larman \cite{larm84} raised the following combinatorial problem:

\medskip
\noindent
{\bf Larman's Problem:} {\it Let $\mathcal{A}$ be a family of subsets of $\{1,2, \ldots ,n\}$ such that every
two members of $\mathcal{A}$ overlap in at least $k$ elements. Can $\mathcal{A}$ be divided into $n$ subfamilies $\mathcal{A}_1$,
$\mathcal{A}_2$,  ..., $\mathcal{A}_n$ such that every two members of $\mathcal{A}_i$ overlap in at least $k+1$ elements}?

\medskip
Given $\mathcal{A}$ as in the statement of Larman's problem, let $\ell (\mathcal{A},n,k)$ denote the smallest number
$m$ for which there exist subfamilies $\mathcal{A}_1$, $\mathcal{A}_2$, ..., $\mathcal{A}_m$ such that
$$\mathcal{A} = \bigcup_{i=1}^m\mathcal{A}_i$$
and every two members of $\mathcal{A}_i$ overlap in at least $k+1$ elements.  Then an affirmative answer to Larman's
problem for the integer $n$ implies that $\ell (\mathcal{A},n,k)\leq n$.

At first glance, it is not easy to notice a connection between Borsuk's problem and Larman's. However, they are closely related.
Assume in what follows that every member of the family $\mathcal{A}$ has cardinality $h$.

Denote by $T_n$ the {\it mapping} from $\mathcal{A}$ to $\mathbb{E}^n$ defined by
$$T_n(A) = (x^1, x^2, ..., x^n)\eqno(3.1)$$
where
$$x^i=\left\{
\begin{array}{ll}
0, & i\not\in A,\\
1, & i\in A.
\end{array}\right.
$$
Note that two members $A$ and $A'$ of $\mathcal{A}$ overlap in exactly $j$ elements if and only if
$$\Vert T_n(A)-T_n(A')\Vert ={\sqrt {2(h - j)}}.\eqno(3.2)$$
Consequently, letting $T_n(\mathcal{A}') = \{T_n(A): A \in \mathcal{A}'\}$ for any subfamily $\mathcal{A}'$ of $\mathcal{A}$,
we see that
$$d(T_n(\mathcal{A}')) \leq \sqrt{2(h - k)}\eqno(3.3)$$
with equality occuring if and only if some two members of $\mathcal{A}'$ overlap in exactly $k$ elements.  But then
$$b(T_n(\mathcal{A})) = \ell (\mathcal{A}, n, k),\eqno(3.4)$$
and so an affirmative answer to Borsuk's problem in dimension $n$ implies that $\ell(\mathcal{A}, n, k) \leq n + 1$.

Clearly now, if for a given $n$ we can find $\mathcal{A}$ and $k$ as above such that
$$\ell (\mathcal{A}, n, k)> n+1,\eqno(3.5)$$
then the answer to both Larman's problem for the integer $n$ and Borsuk's problem for dimension $n$ will be \lq\lq no". This is just the starting point of Kahn and Kalai's work.

\medskip
In combinatorics, the structures of finite sets were comparatively well-studied. In 1981, P. Frankl and R. M. Wilson \cite{fran81}
proved the following two lemmas.

\medskip
\noindent {\bf Lemma 2.} {\it Let $p$ be a prime and $\mathcal{F}$ be a family of $(2p - 1)$-element subsets of $\{1,2, ..., n\}$ such that
$${\rm card} \{F\cap F'\}\neq p - 1$$
for every two distinct members $F,\ F'\in \mathcal{F}$.  Then}
$${\rm card} \{\mathcal{F}\}\le {n\choose {p-1}}.$$

\medskip
\noindent {\bf Lemma 3.} {\it For $p$ a prime, let $m(p)$ be the maximum number of $2p$-element subsets
of $\{1,2,  ...,4p\}$ such that no two of them overlap in $p$ elements.  Then}
$$m(p) \leq \frac{1}{2}{4p \choose p}. $$

\medskip
Based on Larman's reformulation of Borsuk's problem and Frankl and Wilson's lemmas, J. Kahn and G. Kalai were able to deduce Theorem 4.

\medskip
\noindent
{\bf Remark 4.} Afterwards,  Kahn and Kalai's breakthrough was simplified by N. Alon \cite{alon94} and improved by several authors, in particular by Hinrichs and Richter \cite{hinr03} to $n\ge 298$. In 2014, A. Bondarenko \cite{bond14} presented a $65$-dimensional counterexample to Borsuk's conjecture. Soon after, T. Jenrich and A. E. Brouwer \cite{jenr14} discovered a $64$-dimensional one. In fact,  A. Bondarenko \cite{bond14} presented a $65$-dimensional set of $416$ points which cannot be partitioned into $83$ sets of smaller diameter, and T. Jenrich and A. E. Brouwer \cite{jenr14} discovered a $64$-dimensional set of $352$ points that cannot be divided into fewer than $71$ parts of smaller diameter.

\medskip
Up to now, several upper estimates of $b(X)$ depending only on the dimension $n$ of the set $X$ are known.  In fact, all of them were discovered before Kahn and Kalai's counterexamples. In 1961, L. Danzer \cite{danz61} showed that
$$b(X)<{\sqrt {(n+2)^3(2+{\sqrt 2})^{n-1}\over 3}}.\eqno(3.6)$$
In 1982, M. Lassak \cite{lass82} proved that
$$b(X)\le 2^{n-1}+1.\eqno(3.7)$$
In 1988, by considering sets of constant width, O. Schramm \cite{schr88} (also see Bourgain and Lindenstrass \cite{bour91}) was able to improve these
upper bounds to
$$b(X)\le 5n^{3\over 2}(4+\log n)\left({3\over 2}\right)^{n\over 2}.\eqno(3.8)$$

\vspace{1cm}
\noindent
{\LARGE\bf 4. Borsuk's Problem in Metric Spaces}

\bigskip
\noindent
Let $\mathbb{R}^n$ be an $n$-dimensional linear space over real numbers and let $\sigma$ be a metric defined on $\mathbb{R}^n$. Then $\mathbb{M}^n=\{ \mathbb{R}^n, \sigma\}$, the space together with the metric, is an $n$-dimensional metric space. It is natural to consider Borsuk's problem in general metric spaces.

It is well-known in metric geometry that, if $\sigma$ is a metric defined on $\mathbb{R}^n$, then the set
$$C=\left\{ {\bf x}: {\bf x}\in \mathbb{R}^n,\ \sigma ({\bf o}, {\bf x})\le 1\right\}\eqno(4.1)$$
is a centrally symmetric convex body centered at the origin, usually known as the unit domain of the metric space. On the other hand, if $C$ is a centrally symmetric convex body centered at ${\bf o}$ and ${\bf x}$ and ${\bf y}$ are two points of $\mathbb{R}^n$, defining $\sigma ({\bf x}, {\bf y})$ to be the smallest positive number $\rho $ such that ${1\over \rho }({\bf x}-{\bf y})\in C$, one can easily verify that $\sigma ({\bf x}, {\bf y})$ is a metric defined on $\mathbb{R}^n$. Therefore, in $\mathbb{R}^n$, there is an one-to-one correspondence between metrics and the centrally symmetric convex bodies centered at the origin. For example, the Euclidean metric corresponding to the unit ball, the $\ell_1$ metric corresponding to a cross polytope, and the $\ell_\infty$ metric corresponding to a unit cube.

In 1957, according to Gr\"unbaum \cite{grun57}, it was proved by E. Shamir that, {\it if $\mathbb{M}^2$ is a metric plane such that its unit domain is not a parallelogram, then every bounded set can be separated into three subsets of smaller diameter; if $\mathbb{M}^2$ is a metric plane such that its unit domain is a parallelogram, then every bounded set can be separated into four subsets of smaller diameter}.

In 1957, H. Hadwiger\cite{hadw57} made the following related conjecture:

\medskip
\noindent
{\bf Hadwiger's covering conjecture.} {\it In the $n$-dimensional Euclidean space $\mathbb{E}^n$, every convex body $K$ can be covered by $2^n$ translates of $\lambda K$, where $\lambda$ is a positive number satisfying $\lambda <1$.}

\medskip
It is easy to show that $\lambda K$ can be replaced by ${\rm int} (K)$. The conjecture is simple sounding and its two-dimensional case had been proved by F. W. Levi \cite{levi54} before the conjecture was made. In fact, he showed that {\it there is a positive number $\lambda <1$ such that every parallelogram $P$ can be covered by four translates of $\lambda P$ and every other convex domain $K$ can be covered by three translates of $\lambda K$}. Hadwiger's conjecture has been studied by many authors including K. Bezdek, V. G. Boltjanski, I. T. Gohberg, M. Lassak, H. Martini, C. A. Rogers, V. Soltan and C. Zong. Many partial results are known. For example, any $n$-dimensional convex body $K$ with smooth boundary can be covered by $n+1$ translates of $\lambda K$, where $\lambda$ is a suitable positive number satisfying $\lambda <1$. However, up to now, no complete solution is known for any other dimension.

Assume that $X$ is a bounded set in the metric space $\mathbb{M}^n$ with metric $\sigma $ and let $\widehat{X}$ denote its closed convex hull. For convenience, let $b_\sigma (X)$ denote the smallest number $k$ such that $X$ can be divided into $k$ subsets of smaller diameter with respect to $\sigma $ and let $h(\widehat{X})$ denote the smallest number of translates of $\lambda \widehat{X}$ which can cover $\widehat{X}$, where $\lambda $ is any positive number satisfying $\lambda <1$. It is easy to see that
$$b_\sigma (X)\le h(\widehat{X})\eqno(4.2)$$
holds for all metrics and all bounded sets $X$ in $\mathbb{R}^n$. In 1965, V. G. Boltyanski and I. T. Gohberg \cite{bolt85} proposed the following two problems related to Borsuk's conjecture.

\medskip\noindent
{\bf Problem 1.} Is it true that
$$b_\sigma (X)\le 2^n$$
holds for all bounded sets $X$ in $\mathbb{R}^n$ and all metrics $\sigma $ on $\mathbb{R}^n$?

\medskip\noindent
{\bf Problem 2.} Assume that $C$ is the centrally symmetric convex body determined by the metric $\sigma $ in $\mathbb{R}^n$. Is it true that
$$b_\sigma (X)\le h(C)$$
holds for all bounded sets $X$ in $\mathbb{R}^n$?

\medskip
Clearly, the counterexamples to Borsuk's conjecture listed in Section 3 did provide negative answer to Problem 2 in high dimensions. In 2008, C. Zong \cite{zong08} discovered a particular set $X$ and a centrally symmetric convex body $C$ in $\mathbb{R}^3$ satisfying both
$$b_\sigma (X)=5\eqno(4.3)$$
and
$$h(C)=4,\eqno(4.4)$$
which provides a negative answer for Problem 2 in three dimensions.

\medskip
In 2009, L. Yu and C. Zong \cite{yuzo09} studied Problem 1 and obtained the following partial results.

\medskip\noindent
{\bf Theorem 5.} {\it In three-dimensional $\ell_p$ space
$$b_{\ell_p}(X)\le 2^3$$
holds for all bounded sets $X$.}

\medskip\noindent
{\bf Remark 5.} Clearly, Hadwiger's conjecture implies Theorem 5. However, the conjecture is still open in three dimensions. A computer proof programm was proposed by C. Zong \cite{zong10} in 2010. The centrally symmetric case was proved by M. Lassak \cite{lass84} in 1984.

\medskip
Let $C_p$ denote the unit domain of the three-dimensional $\ell_p$ space.  Let $\tau (p)$ denote the smallest number such that there exists a parallelepiped $P$ satisfying
$$P\subseteq C_p\subseteq \tau (p) P.\eqno(4.5)$$
It can be shown that $\tau (p)\le 2$, where the equality holds if and only if $p=1$. Then theorem 5 can be deduced by considering two cases with respect to $p=1$ and $p>1$.

\medskip\noindent
{\bf Theorem 6.} {\it In $n$-dimensional $\ell_p$ spaces, let $C_p$ denote the unit domain, for every bounded centrally symmetric set $X$ we have

$$b_{\ell_p}(X)\le h(C_p) \le \left\{\begin{array}{ll}
2n& \mbox{if $p=1$,}\\
n+1& \mbox{if $1< p<\infty,$}\\
2^n& \mbox{if $p=\infty $.}
\end{array}
\right.$$}

\medskip
Let $C$ be the unit domain of an $n$-dimensional metric space $\mathbb{M}^n=\left\{\mathbb{R}^n,\sigma \right\}$. For every bounded centrally symmetric set X, one can deduce that
$$b_\sigma (X)\le b_\sigma (C)\le h(C).\eqno(4.6)$$
Then theorem 6 can be shown by considering three cases with respect to $p=1$, $1<p<\infty$ and $p=\infty$.

\vspace{1cm}
\noindent
{\LARGE\bf 5. A Reformulation for Borsuk's Problem}

\bigskip
\noindent
Let $B$ denote the $n$-dimensional unit ball centered at the origin of $\mathbb{E}^n$ and let $\mathcal{K}^n$ denote the space of all $n$-dimensional convex bodies associated with the {\it Hausdorff metric} $\delta^H( \cdot )$, where
$$\delta^H(K_1, K_2)=\min\left\{ r:\ K_1\subset K_2+rB,\ K_2\subset K_1+rB \right\}.\eqno(5.1)$$

\medskip\noindent
{\bf Definition 3.} Let $m$ be a fixed positive integer. For an $n$-dimensional convex body $K$ we define $f_m(K)$ to be the smallest positive number $\theta $ such that $K$ can be divided into $m$ subsets $X_1$, $X_2$, $\ldots $, $X_m$ satisfying
$$d(X_i)\le \theta d(K),\quad i=1, 2, \ldots, m.$$

\medskip
Assume that $K_1$ and $K_2$ are $n$-dimensional convex bodies satisfying $$d(K_1)\ge 2,\eqno(5.2)$$ $$d(K_2)\ge 2,\eqno(5.3)$$ and
$$\delta^H(K_1, K_2)\le \epsilon,\eqno(5.4)$$
where $\epsilon $ is a small positive number. Clearly by (5.4) we have
$$d(K_2)-2\epsilon \le d(K_1)\le d(K_2)+2\epsilon\eqno(5.5)$$
and
$$K_2\subseteq K_1+\epsilon B.\eqno(5.6)$$

If $K_1$ can be divided into $m$ subsets $X_1$, $X_2$, $\ldots $, $X_m$ such that
$$d(X_i)\le \theta_1d(K_1)\eqno(5.7)$$
holds for all $i=1,$ $2,$ $\ldots$, $m$, where $\theta_1=f_m(K_1)<1$. Then, by (5.6) we have
$$K_2=\bigcup_{i=1}^m\left(K_2\cap (X_i+\epsilon B)\right)\eqno(5.8)$$
and
$$d\left(K_2\cap (X_i+\epsilon B)\right)\le d(X_i+\epsilon B)\le \theta_1d(K_1)+2\epsilon \le \theta_1(d(K_2)+2\epsilon )+2\epsilon \le (\theta_1+2\epsilon )d(K_2).\eqno(5.9)$$
Consequently, we get
$$\left| f_m(K_1)-f_m(K_2)\right|\le 2\epsilon \eqno(5.10)$$

In conclusion, we have proved the following lemma.

\medskip\noindent
{\bf Lemma 4.} {\it The functional $f_m(K)$ is continuous on $\mathcal{K}^n$. In particular, when $d(K_1)\ge 2$, $d(K_2)\ge 2$ and $\delta^H(K_1,K_2)\le \epsilon,$ we have}
$$\left| f_m(K_1)-f_m(K_2)\right|\le 2\epsilon .$$

\medskip
In 1958, H. G. Eggleston \cite{eggl58} proved the following result.

\medskip\noindent
{\bf Lemma 5.} {\it Assume that $K$ is an $n$-dimensional convex body of unit constant width. First, its insphere $S_1$ and circumsphere $S_2$ are concentric. Let $r$ and $R$ be the radii of $S_1$ and $S_2$, respectively, then we have}
$$1-\sqrt{ n/(2n+2)}\le r\le R\le \sqrt{ n/(2n+2)}.$$

\medskip
By Lemma 2 and Lemma 5, to solve Borsuk's problem in $\mathbb{E}^n$, it is sufficient to deal with all the convex bodies $K$ satisfying
$$B\subseteq K\subseteq r_nB,\eqno(5.11)$$
where
$$r_n={{\sqrt{ n/(2n+2)}}\over {1-\sqrt{ n/(2n+2)}}}. \eqno(5.12)$$
For convenience, we denote the set of all $n$-dimensional convex bodies satisfying this condition by $\mathcal{D}^n$. Clearly it is a compact connected subset of $\mathcal{K}^n$.

\medskip
With this preparation, Borsuk's problem can be reformulated as following.

\medskip\noindent
{\bf Borsuk's Problem.} Is there a positive $\alpha_n<1$ such that
$$f_{n+1}(K)\le \alpha_n$$
holds for all $K\in \mathcal{D}^n$?

\medskip\noindent
{\bf Remark 6.} In any metric space, Borsuk's corresponding problem can be reformulated in a similar way.

\vspace{1cm}
\noindent
{\LARGE\bf 6. A Computer Program for Borsuk's Conjecture}

\bigskip\noindent
{\bf Definition 4.} Let $\beta $ be a given positive number, and let $K_1,$ $K_2$, $\cdots$, $K_{\varpi(\beta )}$ be $\varpi (\beta )$ convex bodies in $\mathcal{D}^n$. If for each $K\in \mathcal{D}^n$ we always can find a corresponding $K_i$ satisfying
$$\delta^H(K, K_i) \le \beta ,$$
we call $\mathcal{N}=\{ K_1, K_2, \cdots, K_{\varpi(\beta )}\}$ a $\beta$-net in $\mathcal{D}^n$.

\medskip\noindent
{\bf Remark 7.} Writing
$$\mathcal{B}(K_i,\beta )=\left\{ K\in \mathcal{K}^n: \ \delta^H(K, K_i)\le \beta \right\},$$ it is easy to show that $\mathcal{N}=\{
K_1, K_2, \cdots , K_{\varpi(\beta )}\}$  is a $\beta$-net in $\mathcal{D}^n$ if and only if
$$\mathcal{D}^n \subseteq \bigcup_{i=1}^{\varpi(\beta )}\mathcal{B}(K_i, \beta ).$$

\medskip
Let $\mathbb{Z}^n$ be the integer lattice in $\mathbb{E}^n$, let $\kappa$ be a small positive number, and let $\mathbb{P}^n$ denote the set of all lattice polytopes of $\kappa\mathbb{Z}^n$ which are elements of $\mathbb{D}^n$. Assume that $K$ is a convex body in $\mathcal{D}^n$ with boundary $\partial (K)$. For each ${\bf x}\in \partial (K)$ we choose $g({\bf x})$ to be one of its nearest lattice points of $\kappa\mathbb{Z}^n$ and define
$$P={\rm conv}\left(\bigcup_{{\bf x}\in\partial (K)}g({\bf x}) \right).\eqno(6.1)$$
By routine argument, it can be shown that
$$\delta^H(K, P)\le \sqrt{n}\kappa.\eqno(6.2)$$
Therefore, all the lattice polytopes of $\kappa\mathbb{Z}^n$ in $\widehat{\mathcal{D}^n}$ form a $\sqrt{n}\kappa$-net in $\mathcal{D}^n$, where $\widehat{\mathcal{D}^n}$ denotes the set of all lattice polytopes $P$ satisfying
$$(1-\sqrt{n}\kappa)B\subseteq P\subseteq (r_n+\sqrt{n}\kappa )B.\eqno(6.3)$$

\begin{figure}[ht]
\centerline{\includegraphics[height=5cm,width=5.5cm,angle=0]{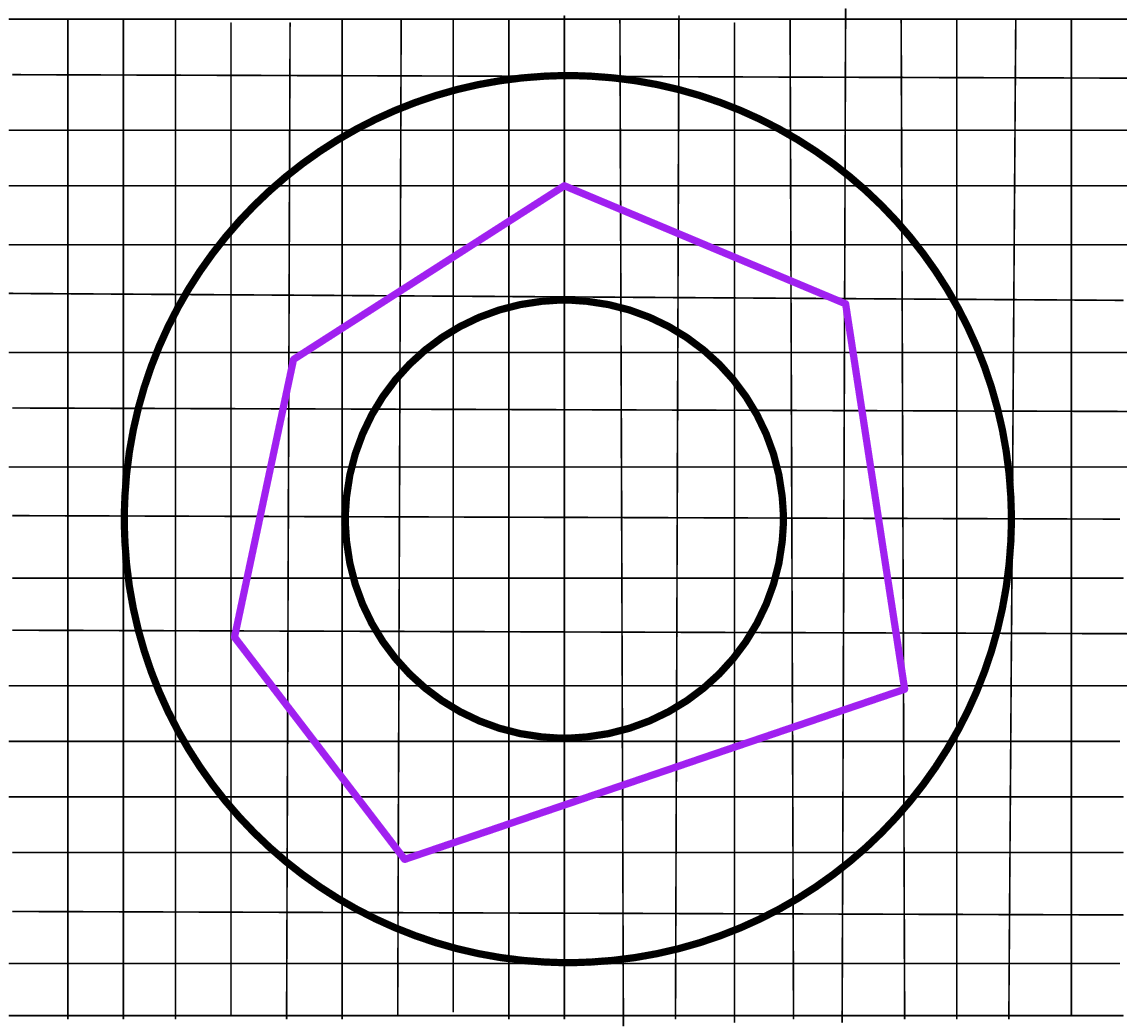}}
\end{figure}

\bigskip\noindent
{\bf A Possible Proof Program}

\medskip\noindent
{\bf Step 1.} Assume that Borsuk's conjecture is true in $\mathbb{E}^n$. Based on some particular examples, one can guess a possible constant $\alpha_n$ such that
$$f_{n+1} (K)\le \alpha_n \eqno(6.4)$$
holds for all $K\in \mathcal{D}^n$.

\medskip\noindent
{\bf Step 2.} Taking
$$\kappa={{1-\alpha_n}\over {4\sqrt{n}}}\eqno(6.5)$$
and defining $\Lambda =\kappa \mathbb{Z}^n$, then for every convex body $K$ in $\mathcal{D}^n$ there is a lattice polytope $P$ in $\widehat{\mathcal{D}^n}$ satisfying
$$\delta^H(P,K)\le {1\over 4}(1-\alpha_n).\eqno(6.6)$$

\medskip\noindent
{\bf Step 3.} Enumerate all the lattice polytopes in $\widehat{\mathcal{D}^n}$. The number of the lattice polytopes is huge. The enumeration can be done only by a computer. For example, by deleting the lattice vertices successively.

\medskip\noindent
{\bf Step 4.} For each lattice polytope $P$, by trying suitable patterns with the help of computer to verify that
$$f_{n+1}(P)\le \alpha_n.\eqno(6.7)$$

\medskip\noindent
{\bf Conclusion.} By Lemma 4, (6.6) and (6.7), one has
$$f_{n+1}(K)\le f_{n+1}(P)+{1\over 2}(1-\alpha_n)\le \alpha_n+{1\over 2}(1-\alpha_n)={1\over 2}(1+\alpha_n)<1.$$
Then, the theorem will follow.

\medskip\noindent
{\bf Example 1.} In $\mathbb{E}^4$, we may try $\alpha_4=0.995$, $r_4=(\sqrt{10}+2)/3$ and $\kappa =0.000625$. Then, we get an huge number (explicitly bounded) of four dimensional lattice polytopes to enumerate and to verify.

\vspace{0.6cm}\noindent
{\bf Acknowledgements.} This work is supported by the National Natural Science Foundation of China (NSFC11921001) and the National Key Research and Development Program of China (2018YFA0704701).

\bibliographystyle{amsplain}

\bigskip\medskip\noindent
Chuanming Zong, Center for Applied Mathematics, Tianjin University, Tianjin, China.

\noindent
Email: cmzong@math.pku.edu.cn

\end{document}